\documentclass[fleqn]{mat01}
\usepackage{times,mathtimy,amssymb,latexsym,graphicx}
\begin{document}

\setcounter{page}{31} \firstpage{31}

\newtheorem{theore}{Theorem}
\renewcommand\thetheore{\arabic{section}.\arabic{theore}}
\newtheorem{definit}[theore]{\rm DEFINITION}
\newtheorem{theor}[theore]{\bf Theorem}
\newtheorem{propo}[theore]{\rm PROPOSITION}
\newtheorem{lem}[theore]{\it Lemma}
\newtheorem{rem}[theore]{\it Remark}

\newtheorem{coro}[theore]{\rm COROLLARY}
\newtheorem{probl}[theore]{\it Problem}
\newtheorem{exampl}[theore]{\it Example}
\newtheorem{pot}[theore]{\it Proof of Theorem}

\def\nota{\trivlist \item[\hskip \labelsep{\it Notations.}]}

\renewcommand\theequation{\thesection\arabic{equation}}






\title{On the stability of Jensen's functional equation on groups}

\markboth{Valeri\u{i} A Fa\u{i}ziev and Prasanna K Sahoo}{Jensen's
functional equation on groups}

\author{VALERI\u I A FA\u IZIEV and PRASANNA K SAHOO$^{*}$}

\address{Tver State Agricultural Academy, Tver Sakharovo, Russia\\
\noindent $^{*}$Department of Mathematics, University of Louisville,
Louisville, KY~40292, USA\\
\noindent E-mail: vfaiz@tvcom.ru; sahoo@louisville.edu}

\volume{117}

\mon{February}

\parts{1}

\pubyear{2007}

\Date{MS received 17 March 2005}

\begin{abstract}
In this paper we establish the stability of Jensen's functional
equation on some classes of groups. We prove that Jensen equation
is stable on noncommutative groups such as metabelian groups and
$T(2, K)$, where $K$ is an arbitrary commutative field with
characteristic different from two. We also prove that any group
$A$ can be embedded into some group $G$ such that the Jensen
functional equation is stable on $G$.
\end{abstract}


\maketitle

\section{Introduction}

Given an operator $T$ and a solution class $\{ u \}$ with the
property that $T(u) = 0$, when does $\|T(v)\|\leq \varepsilon$ for
an $\varepsilon > 0$ imply that $\|u-v\|\leq \delta (\varepsilon)$
for some $u$ and for some $\delta > 0$? This problem is called the
stability of the functional transformation. A great deal of work
has been done in connection with the ordinary and partial
differential equations. If $f$ is a function from a normed vector
space into a Banach space, and $\|f(x+y)-f(x)-f(y)\|\leq
\varepsilon$, Hyers in 1941 proved that there exists an additive
map $A$ such that $\|f(x) - A(x)\|\leq \varepsilon$. If $f(x)$ is
a real continuous function of $x$ over $\mathbb{R}$, and
$|f(x+y)-f(x)-f(y)| \leq \varepsilon$, it was shown by Hyers and
Ulam in 1952 that there exists a constant $k$ such that
$|f(x)-kx|\leq 2 \varepsilon$. Taking these results into account,
we say that the additive Cauchy equation $f(x+y) = f(x)+f(y)$ is
stable in the sense of Hyers and Ulam. For more on stability of
homomorphisms, the interested reader is referred to
\cite{Ulam1,Hyers1,Hyers4,Hyers2,Hyers3} and \cite{Aczel1}. After
Hyers's 1941 result a great number of papers on the subject have
been published, generalizing Ulam's problem and Hyers's theorem in
various directions (see \cite{Forti2},
\cite{Hyers2,Hyers3,Hyers5,Hyers6,Hyers7} and
\cite{Szekelyhidi1}).

In this paper we study the stability of Jensen's functional
equation
\begin{equation*}
f(xy)+f(xy^{-1}) = 2f(x)
\end{equation*}
on some classes of noncommutative groups. This Jensen's equation
was studied in the papers \cite{Aczel2}, \cite{Chung1} and
\cite{Ng2}. The question of stability of this equation was
investigated in \cite{Kominek1,Jung1,Laczkovich1,Tabor1} and
\cite{Yang1}. In all these papers domain of $f$ is either an
abelian group or some of its subsets.

\section{Auxiliary results}

Suppose that $G$ is an arbitrary group and $E$ is an arbitrary
real Banach space. In this sequel, we will write the arbitrary
group $G$ in multiplicative notation so that $1$ will denote the
identity element of $G$.

\begin{definit}$\left.\right.$\vspace{.5pc}

\noindent {\rm We will say that a function $f\hbox{\rm :}\ G\to E$
is a $(G;E)$-{\it Jensen function} if for any $x,y\in G $ we have
\begin{equation}
\label{j0} f(xy)+f(xy^{-1})-2f(x)=0.
\end{equation}
We denote the set of all $(G;E)$-Jensen functions by $J(G;E)$.}
\end{definit}

Denote by $J_0(G;E)$ the subset of $J(G;E)$ consisting of
functions $f$ such that $f(1)=0$. Obviously  $J_0(G;E)$ is a
subspace of $J(G;E)$ and $J(G;E)=J_0(G;E)\oplus E$.

\begin{definit}$\left.\right.$\vspace{.5pc}

\noindent {\rm We will say that a function $f\hbox{\rm :}\ G\to E$
is a {$(G;E)$-quasi-Jensen function} if there is $c>0$ such that
for any $x,y\in G $ we have
\begin{equation}
\label{c} \|f(xy)+f(xy^{-1})-2f(x) \| \le c.
\end{equation}}
\end{definit}

It is clear that the set of ($G;E$)-quasi-Jensen functions is a
linear real space. Denote it by $KJ(G;E)$. From (\ref{c}) we
obtain
\begin{equation*}
\|f(y)+f(y^{-1})-2f(1) \| \le c.
\end{equation*}
Therefore
\begin{equation}
\label{c1} \|f(y)+f(y^{-1})\|  \le c_1 ,
\end{equation}
where $c_1 = c + \|2f(1)\|$. Now letting $x$ for $y$ in  (\ref{c}), we get
\begin{equation*}
\|f(x^2)+f(1)-2f(x) \| \le c.
\end{equation*}
Hence
\begin{equation}
\label{c2} \|f(x^2)-2f(x) \| \le c_2 ,
\end{equation}
where $c_2 = c +\|f(1)\|$. Again substitution of $y=x^2$ in (\ref{c}) yields
\begin{equation*}
\|f(x^3)+f(x^{-1})-2f(x) \|\le c .
\end{equation*}
Thus taking into account~(\ref{c1}) we obtain
\begin{equation}
\label{c3} \|f(x^3)-3f(x) \| \le c_3,
\end{equation}
where $c_3=  c +c_1 $.

Let $c$ be as in~(\ref{c}) and define the set $C$ as follows: $C=\{\,c_m\, | \, \, \,
m\in \mathbb{N} \,\}$, where $c_1=c+2\| f(1)\|$, $c_2=c+\| f(1)\|$, $c_3=c+c_1$ and
$c_m=c+c_1+c_{m-2}$, if $m>3$.

\begin{lem}\label{l-1} Let $f\in KJ(G;E)$ such that
\begin{equation*}
\|\,f(xy)+f(xy^{-1})-2f(x)\,\|\le c.
\end{equation*}
Then for any $x\in G$ and any $m\in \mathbb{N}$ the following relation holds{\rm :}
\begin{equation}
\label{cm} \|f(x^m)-mf(x) \| \le c_m.
\end{equation}
\end{lem}

\begin{proof}
The proof is by induction on $m$. For $m=3$, the Lemma is established. Suppose that
for $m$ the lemma has been already established, let us verify it for $m+1$. Letting
$y=x^m$ in (\ref{c}), we have
\begin{equation*}
\|f(x^{m+1})+f(x^{-(m-1)}-2f(x) \| \le c.
\end{equation*}
From (\ref{c1}) we obtain
\begin{equation*}
\|f(x^{m+1})-f(x^{m-1})-2f(x) \| \le c +c_1 .
\end{equation*}
By induction hypothesis we have
\begin{equation*}
\|f(x^{m-1})-(m-1)f(x) \| \le c_{m-1}
\end{equation*}
and hence,
\begin{equation*}
\|f(x^{m+1})-(m+1)f(x)\| \le c_{m+1}= c +c_1 +c_{m-1}.
\end{equation*}
Now the lemma is proved.\hfill $\Box$
\end{proof}

\begin{lem}\label{l-2}
Let $f \in KJ (G; E)$. For any $m>1${\rm ,} $k\in \mathbb{N}$ and $x\in G$ we have
\begin{equation}\label{dm}
\|f(x^{m^k})-m^kf(x) \| \le c_m(1+m+\cdots +m^{k-1})
\end{equation}
and
\begin{equation}\label{dm1}
\left \|\frac{1}{m^{k}}f(x^{m^{k}})-f(x) \right \|
\le c_m.
\end{equation}
\end{lem}

\begin{proof}
The proof will be based on induction on $k$. If $k=1$, then (\ref{dm}) follows from
(\ref{cm}). Suppose that (\ref{dm}) for $k$ is true, let us verify it for $k+1$.
Substituting $x^m$ for $x$ in (\ref{dm}) implies
\begin{equation*}
\|f(x^{m^{k+1}})-m^kf(x^m) \| \le c_m(1+m+\cdots +m^{k-1}).
\end{equation*}
Now using (\ref{cm}) we obtain
\begin{equation*}
\|m^kf(x^m)-m^{k+1}f(x) \| \le c_mm^{k}
\end{equation*}
and hence
\begin{equation*}
\|f(x^{m^{k+1}})-m^{k+1}f(x) \| \le c_m(1+m+\cdots +m^{k}).
\end{equation*}
The latter implies
\begin{equation*}
\left \|\frac{1}{m^{k+1}}f(x^{m^{k+1}})-f(x) \right \| \le c_m (1+m+\cdots +m^{k})
\frac{1}{m^{k+1}}\le c_m.
\end{equation*}
This completes the proof of the lemma.\hfill$\Box$
\end{proof}

From~(\ref{dm1}) it follows that the set
\begin{equation*}
\left\{\frac{1}{m^{k}}f(x^{m^{k}}) \, | \quad k\in \mathbb{N}
\right\}
\end{equation*}
is bounded.

Substituting $x^{m^n}$ in place of $x$ in~(\ref{dm1}), we obtain
\begin{align*}
\left \|\frac{1}{m^{k}}f(x^{m^{n+k}})-f(x^{m^n}) \right \| &\le
c_m,\\[1pc]
\left \|\frac{1}{m^{n+k}}f(x^{m^{n+k}})-\frac{1}{m^{n}}f(x^{m^n}) \right \| &\le
\frac{c_m}{m^n}\to 0,\quad \hbox{as} \quad n\to \infty .
\end{align*}
From the latter it follows that the sequence
\begin{equation*}
\left\{\frac{1}{m^{k}}f(x^{m^{k}}) \, | \quad k\in \mathbb{N} \right\}
\end{equation*}
is a Cauchy sequence. Since the space $E$ is complete, the above sequence has a limit
and we denote it by $\varphi_m(x)$. Thus
\begin{equation*}
\varphi_m(x)=\lim_{k\to \infty}\frac{1}{m^{k}}f(x^{m^{k}}).
\end{equation*}
From~(\ref{dm1}) it follows that
\begin{equation}
\label{dm2} \|\varphi_m(x)-f(x)\|\le c_m, \quad \forall x\in G.
\end{equation}

\begin{lem}\label{l-3}
Let $f\in KJ(G;E)$ such that
\begin{equation*}
\|f(xy)+f(xy^{-1})-2f(x)\|\le c, \qquad \forall x,y \in G .
\end{equation*}
Then for any $m\in \mathbb{N}$ we have $\varphi_m \in KJ(G;E)$.
\end{lem}

\begin{proof}
Indeed, by~(\ref{dm2})
\begin{align*}
&\|\varphi_m(xy)+\varphi_m(xy^{-1})-2\varphi_m(x) \|\\[.6pc]
&\quad\,=\|\varphi_m(xy)-f(xy) +\varphi_m(xy^{-1}) -f(xy^{-1})-
2\varphi_m(x) +2f(x)\\[.6pc]
&\qquad\,+f(xy)+f(xy^{-1})-2f(x) \|\\[.6pc]
&\quad\,\le \|\varphi_m(xy)-f(xy)\| +\|\varphi_m(xy^{-1}) -f(xy^{-1})\|\\[.6pc]
&\qquad\,+2\|\varphi_m(x) -f(x) \|+\|f(xy)+f(xy^{-1})-2f(x) \|\\[.6pc]
&\quad\,\le 4c_m +c.
\end{align*}
This completes the proof of the lemma.\hfill$\Box$
\end{proof}

For any $x\in G$ we have the relation
\begin{equation}\label{j11}
\varphi_m(x^{m^k})=m^k\varphi_m(x).
\end{equation}
Indeed
\begin{align*}
\varphi_m(x^{m^k})&=\displaystyle{\lim_{\ell\to
\infty}}\frac{1}{m^\ell}f((x^{m^k})^{m^\ell})\\[.7pc]
&=\displaystyle{\lim_{\ell\to
\infty}}\frac{m^k}{m^{k+\ell}}f(x^{m^{k+\ell}})\\[.7pc]
&=m^k \displaystyle{\lim_{p\to
\infty}}\frac{1}{m^{p}}f(x^{m^{p}})\\[.7pc]
&=m^k\varphi_m(x) .
\end{align*}

\begin{lem}\label{l-4}
If $f\in KJ(G;E)${\rm ,} then $\varphi_2=\varphi_m$ for any $m\ge 2$.
\end{lem}

\begin{proof}
By Lemma~\ref{l-3} we have $\varphi_2,\varphi_m \in KJ(G;E)$. Hence the function
\begin{equation*}
g(x)=\lim_{k\to \infty}\frac{1}{m^k}\varphi_2(x^{m^k})
\end{equation*}
is well-defined and is a $(G;E)$-quasi-Jensen function.

It is clear that $g(x^{m^k})=m^kg(x)$ and $g(x^{2^k})=2^kg(x)$ for any $x\in G$ and
any $k\in \mathbb{N}$. From~(\ref{dm2}) it follows that there are $d_1, d_2\in
\mathbb{R}_+ $ such that for all $x\in G$,
\begin{equation}\label{j12}
\|\varphi_2(x) -g(x)\|\le d_1\quad {\rm and} \quad \|\varphi_m(x) -g(x)\|\le d_2.
\end{equation}
Hence $g\equiv\varphi_2$ and $g\equiv\varphi_m$ and we obtain
$\varphi_2 \equiv\varphi_m$. \hfill$\Box$
\end{proof}

\begin{definit}$\left.\right.$\vspace{.5pc}

\noindent {\rm By $(G;E)$-{pseudo-Jensen function\/} we will mean
a $(G;E)$-quasi-Jensen function $f$ such that $f(x^n)=nf(x)$ for
any $x\in G$ and any $n\in \mathbb{Z}$.}
\end{definit}

\begin{lem}\label{l-5}
For any $f\in KJ(G;E)$ the function
\begin{equation*}
\hat{f}(x)=\lim_{k\to \infty}\frac{1}{2^k}f(x^{2^k})
\end{equation*}
is well-defined and is a $(G;E)$-pseudo-Jensen function such that for any $x\in G${\rm
,}
\begin{equation*}
\|\hat{f}(x) -f(x)\|\le c_2.
\end{equation*}
\end{lem}

\begin{proof}
By Lemma \ref{l-3}, $\hat{f}$ is a $(G;E)$-quasi-Jensen function.
Now by Lemma~\ref{l-4}, we\break have $\hat{f}(x^m) =
\varphi_m(x^m)=m\varphi_m(x)= m\hat{f}(x)$. Thus $\varphi_m(x)=
\hat{f}(x)$ and hence $\varphi_2(x)= \hat{f}(x)$ by
Lemma~\ref{l-4}. From equality $\hat{f}=\varphi_2$ we have
$\|\hat{f}(x) -f(x)\|=\|\varphi_2(x) -f(x)\|\le c_2$. \hfill$\Box$
\end{proof}

\begin{rem}\label{rem1}
{\rm If $f\in PJ(G;E)$, then
\begin{enumerate}
\renewcommand\labelenumi{\arabic{enumi}.}
\leftskip -.2pc

\item $f(x^{-n})=-nf(x)$ for any $x\in G$ and $n\in \mathbb{N}$;

\item if $y\in G$ is an element of finite order then $f(y)=0$;

\item if $f$ is a bounded function on $G$, then $f\equiv 0$.
\end{enumerate}}\vspace{-.5pc}
\end{rem}

\begin{proof}
For some $c>0$ the following relation holds:
\begin{equation*}
\|f(xy)+f(xy^{-1}) -2f(x)\|\le c.
\end{equation*}
From (\ref{c1}) it follows that
\begin{equation*}
\|f(y^k)+f(y^{-k})\|\le c_1 \quad \forall y\in G, \, \forall k\in \mathbb{N}.
\end{equation*}
The last inequality is equivalent to $k\|f(y)+f(y^{-1})\|\le c_1$ or
$\|f(y)+f(y^{-1})\|\le \frac{c_1}{k}$ for all $y\in G$ and all $ k\in \mathbb{N}$. The
latter implies $f(y^{-1})=-f(y)$. Thus for any $n\in \mathbb{N}$ we have
$f(y^{-n})=f({(y^n)}^{-1})=-f(y^n)=-nf(y)$. Hence, the assertion {\rm 1} is
established.

Similarly we verify assertions {\rm 2} and {\rm
3}.\hfill$\Box$
\end{proof}

We denote by $B(G;E)$ the space of all bounded functions on a group $G$ that take
values in $E$.

\begin{theor}[\!]\label{directsum}
For an arbitrary group $G$ the following decomposition holds{\rm :}
\begin{equation*}
KJ(G;E)=PJ(G;E)\oplus B(G;E).
\end{equation*}
\end{theor}

\begin{proof}
It is clear that $PJ(G;E)$ and $B(G;E)$ are subspaces of
$KJ(G;E)$, and $PJ(G;E)\cap B(G;E)=\{0\}$. Hence the subspace of
$KJ(G;E)$ generated by $PJ(G;E)$ and $B(G;E)$ is their direct sum.
That is $PJ(G;E) \oplus B(G;E) \subseteq KJ(G;E)$. Let us verify
that $KJ(G;E)\subseteq PJ(G;E) \oplus B(G;E)$. Indeed, if $f\in
KJ(G;E)$, then by Lemma~\ref{l-5} we have $\hat{f}\in PJ(G;E)$ and
$\hat{f}-f \in B(G;E)$.\hfill$\Box$
\end{proof}

\begin{definit}$\left.\right.$\vspace{.5pc}

\noindent {\rm Let $E$ be a Banach space and $G$ be a group. A
mapping $f\hbox{\rm :}\ G\to E$ is said to be a $(G;E)$-{\it
quasiadditive mapping} of a group $G$ if set $\{f(xy) - f(x) -f(y)
| x,y \in G  \} $ is bounded.}
\end{definit}

\begin{definit}$\left.\right.$\vspace{.5pc}

\noindent {\rm By a $(G;E)$-{pseudoadditive mapping} of a group
$G$ we mean its $(G;E)$-quasiadditive mapping $f$ that satisfies
$f(x^n)= nf(x)$ for all $x \in G$ and all $n \in \mathbb{Z}$.}
\end{definit}

\begin{definit}$\left.\right.$\vspace{.5pc}

\noindent {\rm A {quasicharacter} of a group $G$ is a real-valued
function $f$ on $G$ such that the set $\{f(xy) - f(x) -f(y)| x,y
\in G  \} $ is bounded.}
\end{definit}

\begin{definit}$\left.\right.$\vspace{.5pc}

\noindent {\rm By a {pseudocharacter} of a group $G$ we mean its
quasicharacter $f$ that satisfies $f(x^n)= nf(x)$ for all $x \in
G$ and all $n \in \mathbb{Z}$.}
\end{definit}

The set of all $(G;E)$-quasiadditive mappings  is a vector space
(with respect to the usual operations of addition of functions and
their multiplication by numbers), which will be denoted by
$KAM(G;E)$. The subspace of $KAM(G;E)$ consisting of
$(G;E)$-pseudoadditive mappings will be denoted by $PAM(G;E)$ and
the subspace consisting of additive mappings from $G$ to $E$ will
be denoted by $\hbox{Hom}(G;E)$. We say that a
$(G;E)$-pseudoadditive mapping $\varphi$ of the group $G$ is {\it
nontrivial\/} if $\varphi \notin \hbox{Hom}(G;E) $.

The space of quasicharacters will be denoted by $KX(G)$, the space
of pseudocharacters will be denoted by $PX(G)$, and the space of
real additive characters on $G$ will be denoted by $X(G)$.

\begin{rem}
{\rm If a group $G$ has nontrivial pseudocharacter, then for any Banach space $E$
there is nontrivial $(G;E)$-pseudoadditive mapping.}
\end{rem}

\begin{proof}
Let $f$ be a nontrivial pseudocharacter of the group $G$ and $e\in
E$ such that $e\ne 0$. Consider a mapping $\varphi\hbox{\rm :}\
G\to E$ such that $\varphi(x)= f(x)\cdot e$. It easy to see that
$\varphi$ is a nontrivial $(G;E)$-additive
mapping.\hfill$\Box$
\end{proof}

In \cite{Faiziev10} and \cite{Faiziev11} some classes of groups having nontrivial
pseudocharacters are considered.

\begin{theor}[\!]\label{xy=yx}
For any group $G$ the following relations hold{\rm :}
\begin{enumerate}
\renewcommand\labelenumi{{\rm \arabic{enumi}.}}
\leftskip -.2pc
\item $KAM(G;E)\subseteq KJ(G;E)${\rm ,} $PAM(G;E)\subseteq PJ(G;E)${\rm
,} $\hbox{\rm Hom}(G;E)\subseteq J_0(G;E)$.

\item If $f\in PJ(G;E)$ and for any $x,y \in G$ $f(xy)=f(yx)${\rm ,} then
$f\in PAM(G;E)$.
\end{enumerate}
\end{theor}

\begin{proof}$\left.\right.$

\begin{enumerate}
\renewcommand\labelenumi{\arabic{enumi}.}
\leftskip -.2pc

\item Let $f\in KAM(G;E)$ and $c>0$ such that $\|f(xy)-f(x)-f(y)\|\le c$
for all $x,y\in G$. Then we have
\begin{align*}
&\|f(xy)+f(xy^{-1})-2f(x)\|\\[.4pc]
&\quad\,=\|f(xy)-f(x)-f(y)+ f(xy^{-1})-f(x)-f(y^{-1})\\[.4pc]
&\qquad\,+2f(x)+f(y)+f(y^{-1})- 2f(x) \|\\[.4pc]
&\quad\,=\|f(xy)-f(x)-f(y)+f(xy^{-1})-f(x)-f(y^{-1})\|\\[.4pc]
&\quad\,\le \|f(xy)-f(x)-f(y)\| + \|f(xy^{-1})-f(x)-f(y^{-1})\|\\[.4pc]
&\quad\,\le 2c ,
\end{align*}
that is, $KAM(G;E)\subseteq KJ(G;E)$. Hence, $PAM(G;E)\subseteq PJ(G;E)$.

\item Let $f\in PJ(G;E)$, $c>0$ such that
$\|f(xy)+f(xy^{-1})-2f(x)\|\le c$ and $f(xy)=f(yx)$ for all $x,y\in G$. Then we have
\begin{align*}
&2\|f(xy)-f(x)-f(y)\|\\[.4pc]
&\quad\,=\|f(xy)+f(xy^{-1})-2f(x)+ f(xy)+ f(yx^{-1})- 2f(y)\|\\[.4pc]
&\quad\,\le \|f(xy)+f(xy^{-1})-2f(x)\|\\[.4pc]
&\qquad\,+ \|f(xy)+ f(yx^{-1})- 2f(y)\| \le 2c .
\end{align*}
Hence $\|f(xy)-f(x)-f(y)\|\le c$ and $f\in PAM(G;E)$.\hfill$\Box$
\end{enumerate}
\end{proof}

$\left.\right.$\vspace{-3.5pc}

\begin{coro}\label{abelian}$\left.\right.$\vspace{.5pc}

\noindent If $G$ is an abelian group{\rm ,} then $PJ(G;E) = {\rm Hom}(G;E)$.
\end{coro}

\begin{proof}
By Theorem~\ref{xy=yx} we have $PJ(G;E)=PAM(G;E)$. Let $f\in PAM(G;E)$. Then for some
$c>0$ and for any $n\in \mathbb{N}$, and $a,b \in G$ we have
\begin{align*}
n \| f(ab)-f(a)- f(b) \| &= \|\, f((ab)^n) -f(a^n) -f(b^n)\, \|\\[.2pc]
&=\| f(a^nb^n) -f(a^n) -f(b^n)\|\\[.2pc]
&\le c.
\end{align*}
The latter is possible only if $f\in {\rm Hom}(G;E)$.
\hfill$\Box$

\end{proof}

\section{Stability}

Suppose that $G$ is a group and $E$ is a real Banach space.

\setcounter{theore}{0}
\begin{definit}$\left.\right.$\vspace{.5pc}

\noindent {\rm We shall say that eq.~(\ref{j0}) {is stable} for
the pair $(G;E)$ if for any $f\hbox{\rm :}\ G\to E$ satisfying
functional inequality
\begin{equation*}
\|f(xy)+f(xy^{-1})-2f(x)\|\le c \quad \forall x,y \in G
\end{equation*}
for some $c> 0$ there is a solution $j$ of the functional
equation~\eqref{j0} such that the function $j(x)-f(x)$ belongs to
$B(G;E)$.}
\end{definit}

It is clear that eq.~(\ref{j0}) is stable on $G$ if and only if $PJ(G;E)=J_0(G;E)$.
From Corollary~\ref{abelian} it follows that eq.~(\ref{j0}) is stable on any abelian
group. We will say that a $(G;E)$-pseudo-Jensen function $f$ in nontrivial if $f\notin
J_0(G;E)$.

\begin{theor}[\!]\label{st}
Let $E_1${\rm ,} $E_2$ be a Banach space over reals. Then
eq.~$(\ref{j0})$ is stable for the pair $(G;E_1)$ if and only if
it is stable for the pair $(G;E_2)$.
\end{theor}

\begin{proof}
Let $E$ be a Banach space and $\mathbb{R}$ be the set of reals.
Suppose that eq.~(\ref{j0}) is stable for the pair $(G;E)$.
Suppose that (\ref{j0}) is not stable for the pair $(G,
\mathbb{R})$, then there is a nontrivial real-valued pseudo-Jensen
function $f$ on $G$. Now let $e\in E$ and $\|e\|=1$. Consider the
function $\varphi\hbox{\rm :}\ G\to E $ given by the formula
$\varphi(x)=f(x)\cdot e$. It is clear that $\varphi $ is a
nontrivial pseudo-Jensen $E$-valued function, and we obtain a
contradiction.

Now suppose that eq.~$(\ref{j0})$ is stable for the pair
$(G,\mathbb{R})$, that is, $PJ(G,\mathbb{R})=J_0(G,\mathbb{R})$.
Denote by $E^*$ the space of linear bounded functionals on $E$
endowed by functional norm topology. It is clear that for any
$\psi \in PJ(G,H)$ and any $\lambda \in H^*$ the function $\lambda
\circ \psi $ belongs to the space $PJ(G,\mathbb{R})$. Indeed, let
for some $c>0$ and any $x,y\in G$ we have
$\|\psi(xy)+\psi(xy^{-1})-2\psi(x)\|\le c. $ Hence
\begin{align*}
&|\lambda \circ \psi(xy)+\lambda \circ \psi(xy^{-1})- \lambda
\circ \psi(2x) |\\[.2pc]
&\quad\,=| \lambda(\psi (xy) +\psi (xy^{-1})- 2\psi(x))| \le c \|\lambda\|.
\end{align*}
Obviously, $\lambda \circ \psi(x^n)=n\lambda \circ \psi(x)$ for
any $x\in G$ and for any $n \in \mathbb{N}$. Hence the function
$\lambda \circ \psi$ belongs to the space $PJ(G,\mathbb{R})$. Let
$f\hbox{\rm :}\ G\to H$ be a nontrivial pseudo-Jensen mapping.
Then $x,y\in G$ such that $f(xy)+f(xy^{-1})-2f(x) \ne 0$.
Hahn--Banach theorem implies that there is a $\ell \in H^*$ such
that $\ell(f(xy)+f(xy^{-1})-2f(x)) \ne 0 $, and we see that $\ell
\circ f$ is a nontrivial pseudo-Jensen real-valued function on
$G$. This contradiction proves the
theorem.\hfill$\Box$
\end{proof}

In what follows the space $KJ(G, \mathbb{R})$ will be denoted by $KJ(G )$, the space
$PJ(G, \mathbb{R})$ will be denoted by $PJ(G)$, the space $J(G, \mathbb{R} )$ will be
denoted by $J(G,\mathbb{R})$, and the space $J_0(G,\mathbb{R})$ will be denoted by
$J_0(G)$.

\begin{coro}$\left.\right.$\vspace{.5pc}

\noindent Equation~$(\ref{j0})$ over a group $G$ is stable if and only if
$PJ(G)=J_0(G)$.
\end{coro}

Due to the previous theorem we may simply say that eq.~(\ref{j0}) is stable or not
stable.

\begin{rem}\label{rem1a}
{\rm For any group $G$ and any Banach space $E$ the following relation $PAM(G;E)\cap
J(G;E)={\rm Hom}(G;E)$ holds.}
\end{rem}

\begin{proof}
It is clear that ${\rm Hom}(G;E) \subseteq PAM(G;E)\cap J(G;E)$.

Suppose that $f\in PAM(G;E)\cap J(G;E)$. Then by Lemma 1 from \cite{Faiziev9} we have
$f(xy)=f(yx)$. Since $f \in J(G;E)$, the map $f$ satisfies \setcounter{equation}{0}
\begin{equation}\label{rem1-1}
f(xy)+f(xy^{-1})-2f(x)=0.
\end{equation}
Interchanging $x$ with $y$ in \eqref{rem1-1}, we have
\begin{equation*}
f(yx)+f(yx^{-1})-2f(y)=0
\end{equation*}
which is
\begin{equation}\label{rem1-2}
f(xy)-f(xy^{-1})-2f(y)=0.
\end{equation}
Adding (\ref{rem1-1}) and (\ref{rem1-2}) we obtain $2f(xy)-2f(x)-2f(y)=0.$ Hence
$f(xy)=f(x)+f(y)$ and $f\in {\rm Hom}(G;E)$. Thus we obtain
\begin{equation}\label{hom}
PAM(G;E)\cap J(G;E)={\rm Hom}(G;E)
\end{equation}
and the proof is complete.\hfill$\Box$
\end{proof}

\begin{rem}\label{rem1b}
{\rm If a group $G$ has nontrivial pseudocharacter, then eq.~$(\ref{j0})$ is not
stable on $G$.}
\end{rem}

\begin{proof}
Let let $\varphi$ be a nontrivial pseudocharacter of $G$. Suppose
that there is $j\in J_0(G)$ such that the function $\varphi -j$ is
bounded. Then there is $c>0$ such that $|\varphi(x) -j(x)|\le c $
for any $x\in G$. Hence for any $n\in \mathbb{N} $ we have $c\ge
|\varphi(x^n) -j(x^n)|$ $=n|\varphi(x) -j(x)|$ and we see that the
latter is possible if $\varphi(x) =j(x)$. So, $\varphi\in
PX(G)\cap J_0(G)$. Hence, $f\in X(G)$ and we come to a
contradiction with the assumption about
$f$.\hfill$\Box$
\end{proof}

Let $G$ be an arbitrary group. For $a,b,c\in G$, we set $[a,b]=a^{-1}b^{-1}ab$ and
$[a,b,c]=[[a,b],c]$.

\begin{definit}$\left.\right.$\vspace{.5pc}

\noindent {\rm We shall say that $G$ is {\it metabelian} if for any $x,y,z\in G$ we
have $[[x,y],z]=1$.}
\end{definit}

It is clear that if $[x,y]=1$ then $[[x,y],z]=1$, and hence any abelian group is
metabelian.

Our next goal is to prove a stability theorem for any metabelian group.
Consider the group $H$ over two generators $a,b$ and the following
defining relations:
\begin{equation*}
[b,a]a=a[b,a],\quad b[b,a]=[b,a]b.
\end{equation*}
If we set $c=[b,a]$ we get the following representation of $H$
in terms of generators and defining relations:
\begin{equation}\label{groupH}
H= \langle a,b, \| c= [b,a],\quad [c,a]=[c,b]=1 \rangle.
\end{equation}
It is well-known that each element of $H$ can be uniquely represented as
$g=a^mb^nc^k$, where $m,n,k \in {\mathbb{Z}}$. The mapping
\begin{equation*}
g=a^mb^nc^k \rightarrow \begin{bmatrix}
1 & n & k\\
0 & 1 & m\\
0 & 0 & 1
\end{bmatrix}
\end{equation*}
is an isomorphism between $H$ and $UT(3, \mathbb{Z})$.

\begin{lem}\label{metab-1}
Let $f\in PJ(H)$ and $f(c)=0$. Then $f\in PX(H)=X(H)$.
\end{lem}

\begin{proof}
Let $x=a^mb^nc^k$ and $y=a^{m_1}b^{n_1}c^{k_1}$ be two elements from $H$. Then from
the representation~(\ref{groupH}) it follows that
\begin{equation*}
xy=a^{m+m_1}b^{n+n_1}c^{m_1n+ k+k_1}\, ,\quad
yx=a^{m+m_1}b^{n+n_1}c^{mn_1+ k+k_1}.
\end{equation*}
Hence
\begin{align*}
f(xy)&=f(a^{m+m_1}b^{n+n_1})+f(c^{m_1n+ k+k_1})=f(a^{m+m_1}b^{n+n_1}),\\[.2pc]
f(yx)&=f(a^{m+m_1}b^{n+n_1})+f(c^{mn_1+ k+k_1})=f(a^{m+m_1}b^{n+n_1}).
\end{align*}
Thus $f(xy)=f(yx)$ for any $x,y\in G$. By Theorem~\ref{xy=yx} we
obtain that $f\in PX(G)$. From the representation~(\ref{groupH})
it follows that the subgroup of $H$ generated by element $c$ is
the commutator subgroup of $H$. Lemma 2 from~\cite{Faiziev9}
establishes that if $\varphi\in PX(G)$ such that
$f\big|_{G^\prime} \equiv 0$, then $\varphi\in X(G)$. Hence, $f\in
X(H)$ and $PX(G)=X(G)$. \hfill$\Box$\vspace{1pc}
\end{proof}

\begin{lem}\label{metab-2}
Let $f\in PJ(H)$ and $f(a)=f(b)=f(c)=0$. Then $f\equiv 0$.
\end{lem}

\begin{proof} By Lemma~{\ref{metab-1}} we have
$f(a^mb^nc^k)=f(a^m)+f(b^n)+f(c^k) =0$. \hfill$\Box$\vspace{1pc}
\end{proof}

\begin{lem}\label{metab-3}
A function $\phi$ defined by the formula $\phi(a^mb^nc^k)=mn-2k$ is an element of
$J_0(G)$.
\end{lem}

\begin{proof}
It is clear that $\phi(1)=0$. Now let $x=a^mb^nc^k$, $y=a^{m_1}b^{n_1}c^{k_1}$, then
$xy^{-1}=a^mb^nc^k c^{-k_1}b^{-n_1}a^{-m_1}$ $=a^{m-m_1}b^{n-n_1}c^{m_1n_1-m_1n+
k-k_1}$. Hence
\begin{align*}
&f(xy)+f(xy^{-1})-2f(x)\\[.4pc]
&\quad\,=(m+m_1)(n+n_1)-2(m_1n+k+k_1) +(m-m_1)(n-n_1)\\[.4pc]
&\qquad\,-2(m_1n_1 - m_1n +k-k_1) -2(mn-2k)\\[.4pc]
&\quad\, =0
\end{align*}
and the proof of the lemma is now complete.
\hfill$\Box$\vspace{1.5pc}
\end{proof}

\begin{lem}
$PJ(H)=J_0(H)$.\label{metab-4}
\end{lem}

\begin{proof}
Let $g\in PJ(H)$ and $g(a)=\alpha$, $g(b)=\beta$, $g(c)=\gamma$.
Then there are $\psi\in X(H)$ and $\lambda \in \mathbb{R}$ such
that $\psi(a)=\alpha$, $\psi(b)=\beta$, and $\lambda \phi(c)=
\gamma$. Furthermore, we have $j= \psi+\lambda \phi\in J_0(H)$ and
$(g-j)(a)=(g-j)(b)=(g-j)(c)=0$. By Lemma~\ref{metab-2} we get
$(g-j)\equiv 0$. Hence $g=j$ and $g\in J_0(H)$.
\hfill$\Box$
\end{proof}

\begin{theor}[\!] 
Equation $(\ref{j0})$ is stable on any metabelian group.
\end{theor}

\begin{proof}
Let $G$ be a metabelian group and $f\in PJ(G)$. Let $x,y \in G$.
Then there is a homomorphism $\tau$ of $H$ into $G$ such that
$\tau(a)=x $ and $\tau(b)=y$. Obviously, the function
$f^*(g)=f(\tau(g))$ belongs to $PJ(H)$. Now if
$f(xy)+f(xy^{-1})-2f(x)\ne 0$, then
$f^*(ab)+f^*(ab^{-1})-2f^*(a)\ne 0$ and we arrive at a
contradiction with the previous lemma. Thus $f\in J_0(G)$,
$PJ(G)=J_0(G)$ and the eq.~(\ref{j0}) is stable on $G$.
\hfill$\Box$
\end{proof}


\begin{definit}$\left.\right.$\vspace{.5pc}

\noindent {\rm Let $G$ be a group, $f\in PJ(G;E)$, and $b$ an automorphism of $G$. We
will say that $f$ is invariant relative to $b$ if for any $x\in G$ the relation
$f(x^b)=f(x)$ holds. If the latter relation is valid for any $b \in B$, where $B$ is a
group of automorphism of $G$, then we will say that $f$ is invariant relative to $B$.}
\end{definit}

\begin{lem}\label{emb-1}
Let $f$ be an element from $PJ(G;E)$ and $b$ an element of
order two from $G$. Then $f$ is invariant relative to
inner automorphism of $G$ corresponding to element $b$.
\end{lem}

\begin{proof}
Let $\|f(xy)+f(xy^{-1})-2f(x)\|\le c$ for some $c>0$ and
for any $x,y\in G$. Then we have
\begin{align*}
&\|f(bxb)+f(bb^{-1}x^{-1})-2f(b)\|\le c\\[.4pc]
&\|f(bxb)+f(x^{-1})-2f(b)\|\le c\\[.4pc]
&\|f(x^b)+f(x^{-1})\|\le c\\[.4pc]
&\|f(x^b)-f(x)\|\le c.
\end{align*}
From the latter we obtain $\|f({x^n}^b)-f(x^n)\|\le c$ for any $n
\in \mathbb{N}$. Therefore  $n\|f({x}^b)-f(x)\|\le c$ and we get
$f(x^b)=f(x)$. \hfill$\Box$
\end{proof}

Let $K$ be an arbitrary commutative field. Let $K^*$ be the set of nonzero elements of
$K$ with operation of multiplication. Denote by $G$ the group $T(2,K)$ consisting of
matrices of the form
\begin{equation*}
\begin{bmatrix}
\alpha &t\\
0 &\beta\end{bmatrix}; \quad  \alpha, \beta \in K^*; \quad t\in K.
\end{equation*}
Denote by $T$, $E$, $D$ the subgroups of $G=T(2,K)$ consisting of matrices
\begin{equation*}
\begin{bmatrix}
1 &t\\
0 &1
\end{bmatrix}, \quad \begin{bmatrix}
\pm 1 &0\\
0 &\pm 1
\end{bmatrix},
\quad \begin{bmatrix}
 a &0\\
0 &b
\end{bmatrix},
\end{equation*}
respectively, where $a, b \in K^*$ and $ t\in K$.

It is clear that $T\triangleleft G $ and we have the following semidirect products,
$G=D\cdot T$. Subgroup $C$ of $G$ generated by $T$ and $E$ is a semidirect product
$C=E\cdot T$. Now we prove a stability theorem on the noncommutative group $T(2,K)$.

\begin{theor}[\!]\label{stJ-2}
Let $K$ be an arbitrary commutative field. If the characteristic of $K$ is not equal
to two, then the Jensen functional equation is stable on $G$.
\end{theor}

\begin{proof}
Let $f\in PJ(G)$. Every element of $E$ has order two. Hence, by Lemma~\ref{emb-1} we
have $f^e = f$ for any $e\in E$. Here $f^e$ denotes $f( x^e)$ for $x \in G$, and $x^e$
denotes $e^{-1} x e$. Now from the relation
\begin{equation*}
\begin{bmatrix}
-1 & 0\\
0 & 1
\end{bmatrix} \cdot
\begin{bmatrix}
1 &t\\
0 &1
\end{bmatrix} \cdot \begin{bmatrix}
-1 &0\\
0 &1
\end{bmatrix} = \begin{bmatrix}
1 &-t\\
0 &1 \end{bmatrix},
\end{equation*}
it follows that if $e= \big[ \begin{smallmatrix}
-1 &0\\
0 &1
\end{smallmatrix} \big]$ and
$v = \big[ \begin{smallmatrix}
1 &t\\
0 &1
\end{smallmatrix} \big]$, then
\begin{equation*}
f(v)=f^e(v) =f(v^{-1}) =-f(v).
\end{equation*}
Hence $f\big |_T\equiv 0$. It is clear that the map
\begin{equation*}
\tau\hbox{\rm :}\ \begin{bmatrix}
 \alpha &t\\
0 &\beta
\end{bmatrix} \rightarrow \begin{bmatrix}
 \alpha &0\\
0 &\beta
\end{bmatrix}
\end{equation*}
is a homomorphism of $G$ onto $D$. Let $\varphi = f\big|_D$. Then we can extend
$\varphi $ onto $G$ by the rule $\varphi(g)=\varphi(g^\tau)$. It clear that $\varphi
\in PJ(G)$ and for any $d\in D$ we have the following relation $\varphi(d)=f(d)$.

Now let $\omega(x) = f(x) -\varphi(x)$. So $\omega \in PJ(G) $ and $ \omega_{D\cup
T}\equiv 0.$ Let us show that $\omega \equiv 0$ on $G$. For some nonnegative number
$\delta$, we have
\begin{equation}\label{s-1}
|\omega(xy)+\omega(xy^{-1}) -2\omega(x) | \le \delta
\end{equation}
for any $x, y\in G$. Let $x=au, y=bv$, where $a, b\in D$ and $u,v \in T$. Then
\begin{equation*}
xy=abu^bv,\quad xy^{-1}= ab^{-1}(uv^{-1} )^{a^{-1}}.
\end{equation*}
Now from~(\ref{s-1}) it follows
\begin{equation*}
|\omega(abu^bv)+\omega(ab^{-1}(uv^{-1} )^{a^{-1}}) -2\omega(au) | \le \delta.
\end{equation*}
If we put $a=b$ and $u=1$, then from the last relation we get
\begin{equation*}
|\omega(a^2v)+\omega((v^{-1} )^{a^{-1}}) -2\omega(a) | = |\omega(a^2v) | \le \delta.
\end{equation*}
Taking into account equality $(a^2v)^n =(a^n)^2v^{a^{2(n-1)}}v^{a^{2(n-2)}}\cdots
v^{a^2}v$ we get
\begin{equation*}
n|\omega(a^2v) |=|\omega((a^2v)^n)| =
|\omega((a^n)^2v^{a^{2(n-1)}}v^{a^{2(n-2)}}\cdots v^{a^2}v) | \le \delta,
\end{equation*}
that is
\begin{equation*}
|\omega(a^2v) | \le \frac{1}{n}\delta, \,\, \mbox{for any }\,\, n\in \mathbb{N}.
\end{equation*}
It follows that
\begin{equation}
\omega(a^2v) \equiv 0 \,\, \mbox{for any }\,\, a\in D
\,\, \mbox{and any}\,\, v\in T.
\end{equation}
Now let $x=bu$ be an arbitrary element of $G$. Then $x^2=b^2u^bu$.
Therefore $\omega(x) =\frac{1}{2}\omega(x^2)=0$. So, $\omega
\equiv 0$ on $G$, and $PJ(G)=PJ(D)=J_0(D)$. The proof of the
theorem is now complete. \hfill$\Box$
\end{proof}

\section{The theorem of embedding}

In this setion, we prove that any group $A$ can be embedded into
some group $G$ such that the Jensen functional equation is stable
on $G$. From now on, the set of pseudo-Jensen functions on $G$
invariant relative to $B$ will be denoted by $PJ(G,B;E)$ and if
$E=\mathbb{R}$, then the space $PJ(G,B; \mathbb{R} )$ will be
denoted $PJ(G,B)$.

Let $A$ and $B$ be arbitrary groups. For each $b\in B$ denote by
$A(b)$ a group that is isomorphic to $A$ under  isomorphism $a\to
a(b)$. Denote by $D= A^{(B)}=\prod_{b\in B}A(b)$ the direct
product of groups $A(b)$. It is clear that if
$a_1(b_1)a_2(b_2)\cdots a_k(b_k)$ is an element of $D$, then for
any $b\in B$, the mapping
\begin{equation*}
b^*\hbox{\rm :}\ a_1(b_1)a_2(b_2)\cdots a_k(b_k)\to
a_1(b_1b)a_2(b_2b)\cdots a_k(b_kb)
\end{equation*}
is an automorphism of $D$ and $b\to b^*$ is an embedding of $B$
into $\hbox{Aut}\,D$.

Hence, we can form a semidirect product $G= B\cdot D$. This group
is called {\it the wreath product\/} of the groups $A$ and $B$,
and will be denoted by $G=A\wr B$. We will identify the group $A$
with subgroup $A(1)$ of $D$, where $1\in B$. Hence, we can assume
that $A$ is a subgroup of $D$.

Let us denote, by $C$, the group $\prod_{i\in \mathbb{N}}C_i$,
where $C_i$ is a group of order two with generator $b_i$.

\setcounter{theore}{0}
\begin{theor}[\!]\label{stability}
Let $A$ be an arbitrary group. Then $A$ can be embedded into a
group $G$ such that $PJ(G)=J_0(G)=X(G)$. Hence eq.~$(\ref{j0})$ is
stable on group $G$.
\end{theor}

\begin{proof}
Let $C$ be a group as described above. Let us verify that eq.~$(\ref{j0})$ is stable
on $G=A\wr C$. Denote by $D$ the subgroup of $G$ generated by $A(b),\,\,b\in C$. By
Lemma~\ref{emb-1}  we have that if $f\in PJ(G)$, then $f\big|_D\in PJ(D,C)$. Let $b_i,
\, i=1,2, \dots, k$ be distinct elements from $C$. Then for any $a_i, \, i=1,2, \dots,
k$ the subgroup of $D$ generated by $a_i(b_i), \, i=1,2, \dots, k$ is abelian. Hence
if $u=a_1(b_1)a_2(b_2)\cdots a_k(b_k), \,\, $ $v=\alpha_1(b_1)\alpha_2(b_2)\cdots
\alpha_k(b_k)\in D$ and $f\in PJ(D,C)$, then  by Corollary~\ref{abelian}
\begin{align*}
&|f(uv)+ f(uv^{-1})-2f(u)|\\[.4pc]
&\quad\,=  \left | \displaystyle{\sum_{i=1}^k}
[f(a_i\alpha_i(b_i))+ f(a_i\alpha_i^{-1}(b_i))-2f(a_i(b_i))]  \right |.
\end{align*}
Let $b_i$ for $i\in \mathbb{N}$ be distinct elements from $C$. Let $a,\alpha\in A.$
Consider elements $u_k=a(b_1)a(b_2)\cdots a(b_k)$ and
$v_k=\alpha(b_1)\alpha(b_2)\cdots \alpha(b_k)$. Then by Corollary~\ref{abelian}, for
any $k\in \mathbb{N}$, we have
\begin{align*}
&|f(u_kv_k)+ f(u_kv_k^{-1})-2f(u_k)|\\[.4pc]
&\quad\,= \left|\sum_{i=1}^k [f(a\alpha(b_i))+ f(a\alpha^{-1}(b_i)) -2f(a(b_i))]
\right|.
\end{align*}
By Lemma~\ref{emb-1} we have $f(d(b_i))=f(d(1))$ for any $d\in A$ and for any $i\in
\mathbb{N}$. Let $ r= f(a\alpha(b_i))+ f(a\alpha^{-1}(b_i))-2f(a(b_i))$. Hence
\begin{align*}
&|f(u_kv_k)+ f(u_kv_k^{-1})-2f(u_k)|\\[.4pc]
&\quad\,= \left|\sum_{i=1}^k [f(a\alpha(b_i))+
f(a\alpha^{-1}(b_i)) -2f(a(b_i))] \right|.\\[.4pc]
&\quad\, = |k [f(a\alpha(1))+ f(a\alpha^{-1}(1)) -2f(a(1))]|.\\[.4pc]
&\quad\, = k\cdot| r|.
\end{align*}
Further we have
\begin{equation*}
|f(u_kv_k)+ f(u_kv_k^{-1})-2f(u_k)| \le c .
\end{equation*}
Hence
\begin{equation*}
k|r|\le c ,
\end{equation*}
and
\begin{equation*}
|r| \le  c\,\frac{1}{k}\, \quad \forall k\in \mathbb{N} .
\end{equation*}
The latter is possible only if $r=0$. Thus $f\in J_0(D,B)$. Denote
by $j_b$ the restriction of $f$ to $A(b)$. Let $a$ be an arbitrary
element from $A$. According to the action of $C$ on $D$ we have
\begin{equation*}
f(a(b))= f(a(1)^b)= f^b(a(1))=f(a(1)),
\end{equation*}
that is, $j_b(a(b))= j_1(a(1))$. Hence, there is an element $j$ in $J_0(A)$ such that
$j_b(a(b))=j(a)$ for any $a\in A$ and for any $b\in B$. Therefore, for any
$u=a_1(b_1)a_2(b_2)\cdots a_k(b_k)$ the relation
\begin{equation*}
f(a_1(b_1)a_2(b_2)\cdots a_k(b_k))=\sum_{i=i}^kj(a_i)
\end{equation*}
holds. Let $b,b_1,b_2\in C;\,\,d, d_1,d_2\in D$ and $u=b_1d_1,\,\,v=b_2d_2$. Then
\setcounter{equation}{0}
\begin{align}
\label{jf} &|f(uv)+ f(uv^{-1})-2f(u)| \nonumber\\[.2pc]
&\quad\,= |f(b_1b_2d_1^{b_2}d_2)+
f(b_1b_2^{-1}d_1^{b_2^{-1}}{d_2^{-1}}^{b_2^{-1}})-2f(b_1d_1)|\le c.
\end{align}
Further we have
\begin{align*}
&|f(db)+ f(db^{-1})-2f(d)|\le c,\\[.3pc]
&|f(db)+ f(db)-2f(d)|\le c,\\[.3pc]
&|2f(db)-2f(d)|\le c
\end{align*}
or
\begin{equation*}
|f(db)-f(d)|\le \frac{c}{2}.
\end{equation*}
The latter is equivalent to
\begin{equation}
|f(bd^b)-f(d)|\le \frac{c}{2}.
\end{equation}
Taking into account $f(d^b)=f(d)$ we get
\begin{equation}
\label{jf3} |f(bd)-f(d)|\le \frac{c}{2}.
\end{equation}
Putting $b_1=b_2$  in~(\ref{jf}) we obtain
\begin{equation}
\label{jf1} |f(d_1^{b_2}d_2)+f(d_1^{b_2}{d_2^{-1}}^{b_2})-2f(b_2d_1)|\le c.
\end{equation}
Now taking into account~(\ref{jf3}) and the relation $f(d^b)=f(d)$ we obtain
from~(\ref{jf1})
\begin{align*}
&2c \ge |f(d_1^{b_2}d_2)+f(d_1d_2^{-1})-2f(d_1)\\[.2pc]
&\quad \,= |f(d_1^{b_2}d_2)-f(d_1d_2) +f(d_1d_2)+ f(d_1d_2^{-1})-2f(d_1)|.
\end{align*}
Hence
\begin{align}
\label{jf2} &|f(d_1^bd_2)-f(d_1d_2) |\nonumber\\[.2pc]
&\quad\,=|f(d_1^bd_2)-f(d_1d_2) +[f(d_1d_2) +f(d_1d_2^{-1}) -2f(d_1)]\nonumber\\[.2pc]
&\qquad\,-[f(d_1d_2)+f(d_1d_2^{-1}) -2f(d_1)]|\nonumber\\[.2pc]
&\quad\,\le |f(d_1^bd_2)-f(d_1d_2) +[f(d_1d_2) +f(d_1d_2^{-1}) -2f(d_1)]|\nonumber\\[.2pc]
&\qquad\, +|[f(d_1d_2)+f(d_1d_2^{-1}) -2f(d_1)]|\nonumber\\[.2pc]
&\quad\,\le 2c+c =3c
\end{align}
for any $d_1,d_2\in D$ and any $b\in C$.

Let $b\ne 1$. Then from (\ref{jf2}) we obtain $|f(a_1(b)a_2(1))-f(a_1a_2) |\le 3c$ for
any $a_1,a_2\in A$, that is, $|j(a_1(b))+j(a_2(1))-j(a_1a_2) |\le 3c$ for any
$a_1,a_2\in A$, and $|j(a_1)+j(a_2)-j(a_1a_2) |\le 3c$ for any $a_1,a_2\in A$. Hence,
$j\in PX(A)$. But $PX(A) \cap J_0(A)= X(A)$ and we see that $j\in X(A)$.

Let $\psi=f\big|_D$. Then $\psi$ is an element of $X(D)$ invariant
relative to $C$. Let us extend $\psi $ to $G$ as follows:
$\psi_1(bd)=\psi(d)$. It is easy to see that $\psi_1 \in KX(G)$.
From Theorem~\ref{xy=yx} it follows that $\psi_1 \in KJ(G)$. From
Lemma~\ref{l-5} we see that $\hat{\psi}_1 \in PJ(G)$. Let us
verify that $\hat{\psi}_1(xy)=\hat{\psi}_1(yx)$ for all $x,y\in
G$.

Indeed, by Lemma~\ref{l-5} there is $q>0$ such that
\begin{equation}
\label{1} |\hat{\psi}_1(x) -\psi_1(x)|\le q , \quad \forall x \in
G.
\end{equation}
From the relation $\psi_1\in KX(G)$ we see that for some $\delta>0$,
\begin{equation}
\label{1a} |\psi_1(xy) -\psi_1(x)-\psi_1(y)|\le \delta, \quad \forall x, y \in G.
\end{equation}
This implies that for $x,y,z \in G$ the following relation holds:
\begin{equation}
\label{2} |\psi_1(xyz) -\psi_1(x)-\psi_1(y)-\psi_1(z)|\le 2\delta,
\end{equation}
From~(\ref{1}) and~(\ref{2}) we have
\begin{equation}
\label{3} |\hat{\psi}_1((xy)^{n+1}) -\psi_1(x)
-\psi_1((yx)^n)-\psi_1(y)| \le q+2\delta , \ \ \forall x, y \in G.
\end{equation}
Now applying~(\ref{1a}) we obtain
\begin{equation}
\label{3a} |\hat{\psi}_1((xy)^{n+1}) -\psi_1((yx)^{n+1})| \le
q+3\delta , \quad \forall x, y \in G.
\end{equation}
Similarly, we get
\begin{equation}\label{4}
\!\!\!|\hat{\psi}_1((yx)^{n+1}) -\psi_1(x)
-\psi_1((xy)^n)-\psi_1(y)| \le q+2\delta , \ \ \forall x, y \in G.
\end{equation}
Now again using~(\ref{1}), (\ref{1a}) and~(\ref{2}) we obtain
\begin{equation*}
\label{5} |\hat{\psi}_1((xy)^{n+1}) -\psi_1((yx)^n)-\psi_1(yx)|
\le q+3\delta , \quad \forall x,y \in G.
\end{equation*}
From the equality
\begin{equation*}
\hat{\psi}_1((yx)^{n+1})=\hat{\psi}_1((yx)^{n})+\hat{\psi}_1(yx)
\end{equation*}
and~(\ref{1}) we get
\begin{equation}
\label{6}
|\hat{\psi}_1((yx)^{n+1})-\psi_1((yx)^{n})-\psi_1(yx)|\le 2q.
\end{equation}
From~(\ref{4}) and~(\ref{6}) we obtain that for $p=3q+3\delta $ the following
relations hold:
\begin{align*}
&|\hat{\psi}_1((xy)^{n+1})-\hat{\psi}_1((yx)^{n+1})|\le p \quad
\forall x,y \in G \quad \mbox{and}\quad \forall n \in
\mathbb{N}.\\[.2pc]
&(n+1)|\hat{\psi}_1(xy)-\hat{\psi}_1(yx)|\le p.
\end{align*}
This implies that
\begin{equation*}
|\hat{\psi}_1(xy)-\hat{\psi}_1(yx)|\le \frac{p}{n+1} \quad \forall
x,y \in G \quad \mbox{and}\quad \forall n \in \mathbb{N}.
\end{equation*}
The latter is possible only if $\hat{\psi}_1(xy)\equiv
\hat{\psi}_1(yx)$.

Now by Theorem~\ref{xy=yx} we get that $\hat{\psi}_1\in PX(G)$
such that $\hat{\psi}_1\big|_D=\psi$.

It is clear that $g=f-\hat{\psi}_1 \in PJ(G)$ and $g\big|_{C\cup
D}\equiv 0$. Let us verify that $g\equiv 0$ on $G$. Indeed, for
any $bd\in G$ we have $2g(bd)=g((bd)^2)=g(b^2d^bd)=g(d^bd)=0$.
Hence $g(bd)=0$. So, we see that $f=\hat{\psi}_1$ and $f\in
PX(G)$.

Now let us verify that $f\in X(G)$. To do it we verify that $2f$ is a character of the
group $G$. Indeed,
\begin{align*}
&2f(b_1d_1b_2d_2) -2f(b_1d_1)- 2f(b_2d_2)\\[.4pc]
&\quad\,=f(d_1^{b_2b_1b_2}d_2^{b_1b_2}d_1^{b_1}d_2)- f(d_1^{b_1}d_1) - f(d_2^{b_2}d_2)
\\[.4pc]
&\quad\,= f(d_1^{b_1}) +f(d_2^{b_1b_2}) +f(d_1^{b_1}) +f(d_2)
-f(d_1)\\[.4pc]
&\qquad\, -f(d_1)-f(d_2) -f(d_2)\\[.4pc]
&\quad\,\equiv 0.
\end{align*}
So, $2f\in X(G)$ and we obtain that $f\in X(G)$. Hence, $f\in
J_0(G)$ and the equation~(\ref{j0}) is stable on $G$. This
completes the proof. \hfill $\Box$
\end{proof}


\begin{thebibliography}{99}
\bibitem{Aczel1}
Acz\'el J and Dhombres J, Functional equations in several
variables. Encyclopedia of mathematics and its applications
(Cambridge: Cambridge University Press) (1989)
\bibitem{Aczel2}
Acz\'el J, Jung J K and Ng C T, Symmetric second differences in product form on
groups, in: Topics in Mathematical Analysis (ed.) Th~M~Rassias (1989) pp.~1--22
\bibitem{Chung1}
Chung~J~K, Ebanks~B~R, Ng~C~T and Sahoo~P~K, On a
quadratic-trigonometric functional equation and some applications,
{\it Trans. Am. Math. Soc.} {\bf 347} (1995) 1131--1161
\bibitem{Faiziev9}
Fa\u iziev V A, The stability of the equation $f(xy)-f(x)-f(y)=0$ on groups, {\it Acta
Math.} Univ. Comenianae {\bf 1} (2000) 127--135
\bibitem{Faiziev10}
Fa\u iziev V A, Pseudocharacters on a class of extension of free
groups, {\it New~York J. Math.} {\bf 6} (2000) 135--152
\bibitem{Faiziev11}
Fa\u iziev V A, Description of the spaces of pseudocharacters on a free products of
groups, {\it Math. Ineq. Appl.} {\bf 2} (2000) 269--293
\bibitem{Forti2}
Forti G L, Hyers--Ulam stability of functional equations in several variables, {\it
Aequationes Math.} {\bf 50} (1995) 143--190
\bibitem{Hyers1}
Hyers D H, On the stability of the linear functional equation, {\it Proc. Nat. Acad.
Sci. USA} {\bf 27(2)} (1941) 222--224
\bibitem{Hyers4}
Hyers D H, The stability of homomorphisms and related topics, in:
Global Analysis~--~Analysis on Manifolds (ed.) Th~M~Rassias,
Teubner-Texte zur Math., Leipzig (1983) pp.~140--153
\bibitem{Hyers2}
Hyers D H and Ulam S M, On approximate isometry, {\it Bull. Am.
Math. Soc.} {\bf 51} (1945) 228--292
\bibitem{Hyers3}
Hyers D H and Ulam S M, Approximate isometry on the space of
continuous functions, {\it Ann. Math.} {\bf 48(2)} (1947) 285--289
\bibitem{Hyers5}
Hyers D H and Rassias Th~M, Approximate homomorphisms, {\it Aequationes Math.} {\bf
44} (1992) 125--153
\bibitem{Hyers6}
Hyers D H, Isac G and Rassias Th~M, Topics in nonlinear analysis and applications,
(Singapore, New Jersey, London: World Scientific Publ. Co.) (1997)
\bibitem{Hyers7}
Hyers D H, Isac G and Rassias Th~M, Stability of functional
equations in several variables (Boston/Basel/Berlin:
Birkh\"{a}user) (1998)
\bibitem{Jung1}
Jung S M, Hyers-Ulam-Rassias stability of Jensen's equation and
its application, {\it Proc. Am. Math. Soc.} {\bf 126(11)} (1998)
3137--3143
\bibitem{Kominek1}
Kominek Z, On a local stability of the Jensen functional equation, {\it Demonstratio
Math.} {\bf 22} (1989) 199--507
\bibitem{Laczkovich1}
Laczkovich M, The local stability of convexity, affinity and the Jensen equation, {\it
Aequationes Math.} {\bf 58} (1999) 135--142
\bibitem{Yang1}
Lee Y and Jun K, A Generalization of the Hyers-Ulam-Rassias Stability of Jensen's
equation, {\it J. Math. Anal. Appl.} {\bf 238} (1999) 305--315
\bibitem{Ng2}
Ng C T, Jensen's functional equation on groups, {\it Aequationes Math.} {\bf 39}
(1999) 85--99
\bibitem{Szekelyhidi1}
Sz\'ekelyhidi L, Ulam's problem, Hyers's solution~--~and to where
they led, in: Functional equations and inequalities (ed.)
Th~M~Rassias (Kluwer Academic Publishers) (2000) 259--285
\bibitem{Tabor1}
Tabor J and Tabor J, Local stability of the Cauchy and Jensen equations in function
spaces, {\it Aequationes Math.} {\bf 58} (1999) 296--310
\bibitem{Ulam1}
Ulam S M, A collection of mathematical problems (New York:
Interscience Publ.) (1960)

\end{thebibliography}
\end{document}